\documentclass[10pt,a4paper,twoside,openany]{report}
\usepackage[cp1251]{inputenc} 
\usepackage[ukrainian]{babel} 
\usepackage{amsfonts,amsmath,amssymb,amscd}
\usepackage{bezier,graphics,graphicx}
\usepackage{epsfig,layout,latexsym}
\usepackage{theorem}
\usepackage{cite}
\usepackage{indentfirst}
\usepackage{array}
\usepackage{subfig}

\DeclareCaptionLabelSeparator{dot}{. }
\AtBeginDocument{}
\AtBeginDocument{}

\theorembodyfont{\sl}
\newtheorem{theorem}{Theorem}[section]

\newtheorem{corollary}[theorem]{Corollary}
\newtheorem{remark}[theorem]{Remark}
\newtheorem{example}[theorem]{Example}

\numberwithin{equation}{section}
\numberwithin{theorem}{section}

\makeatletter
\renewenvironment{thebibliography}[1]
 {\section*{\centerline{\rm\textsc{Bibliography}}}%
 \@mkboth{\MakeUppercase\refname}{\MakeUppercase\refname}%
 \list{\@biblabel{\@arabic\c@enumiv}}%
 {\settowidth\labelwidth{\@biblabel{#1}}%
 \leftmargin\labelwidth
 \advance\leftmargin\labelsep
 \@openbib@code
  \usecounter{enumiv}%
  \let\p@enumiv\@empty
  \renewcommand\theenumiv{\@arabic\c@enumiv}}%
  \sloppy
  \clubpenalty4000
  \@clubpenalty \clubpenalty
  \widowpenalty4000%
  \sfcode`\.\@m
  \setlength{\itemsep}{-0.1cm}}
  {\def\@noitemerr
  {\@latex@warning{Empty 'thebibliography' environment}}%
 \endlist}
\renewcommand{\@biblabel}[1]{#1.}

\makeatother
\def\R{\mathbb R}
\def\Eg{\mathcal E}
\newcommand{\suml}{\sum\limits}
\newcommand{\intl}{\int\limits}
\DeclareMathOperator{\sinc}{sinc}
\def\Z{\mathbb Z}
\newcommand{\e}{\mathrm{e}}
\newcommand{\bW}{{\mathbf{W}}}
\DeclareMathOperator{\Arg}{Arg}

\begin{document}

\setcounter{equation}{0}
\setcounter{figure}{0}
\setcounter{table}{0}
\setcounter{footnote}{0}
\setcounter{section}{0}

\begin{center}
\textbf{SINC APPROXIMATION OF ALGEBRAICALLY DECAYING FUNCTIONS
}

\setcounter{footnote}{2}
\footnotetext{\textit{Key words}.
sinc methods, sinc interpolation, algebraically decaying functions,
Lambert-W function, polynomial order of convergence, approximation on real-line.}
\end{center}

\def\headlinetitle{SINC APPROXIMATION OF ALGEBRAICALLY DECAYING FUNCTIONS}


\vspace*{10mm}
\centerline{\textsc {D.\,O.\,SYTNYK}}

\def\headlineauthors{D.\,O.\,SYTNYK}

\setcounter{tocdepth}{0}
\addcontentsline{toc}{abcd}{
\textit{D.\,O.\,Sytnyk}\\
Sinc approximation of algebraically decaying functions}

\vspace{4mm}
\begin{small}
\begin{quote}
\textsc{Abstract.}
An extension of sinc interpolation on $\R$ to the class of algebraically decaying functions is developed in the paper. 
Similarly to the classical sinc interpolation we establish two types of error estimates. 
First covers a wider class of functions with the algebraic order of decay on $\R$. 
The second type of error estimates governs the case when the order of function's decay can be estimated everywhere in the horizontal strip of complex plane around $\R$. 
The numerical examples are provided. 
\end{quote}
\end{small}


%
\section*{Introduction}\label{sec:sinc-nonunif_intro}
We begin by introducing some necessary notation. Let
$$
\sinc{(x)} = \frac{\sin{\pi x}}{\pi x},
$$
\begin{equation}\label{eq:Skh}
S\{k,h\}(x) = \sinc\left (\frac{x}{h} - k\right), \quad h>0, \ k \in \Z.
\end{equation}
By $H^1(D_d)$ in the paper we denote the class of functions $f(x)$ analytic in the horizontal strip $D_d$
\begin{equation}\label{eq:Dd} 
D_d=\left\{z=x+iy\: \quad x \in (-\infty,\infty),\
|y| \leq d\right\},
\end{equation}
and such, that the quantity 
$$
N_1(f,D_d) \equiv \int_{\partial D_d} |f(z)| dz,
$$
is bounded. 
Next, for some given $h>0$ and integer $N>0$ we define a sinc interpolation polynomial as 
\begin{equation}\label{eq:CN}
C_N\{f, h\} (x) = \suml_{k=-N}^{N} f(kh)S\{k,h\}(x).
\end{equation}
The following classical result characterize the accuracy of interpolation of $f \in H^1(D_d)$ by $C_N\{f, h\} (x)$  for the case, when $f(s)$ is exponentially decaying.

{\bf Theorem} (Stenger  \cite[p. 137]{Stenger:1993:NMB})
{\it
	Assume that the function  $f \in H^1(D_d)$ is bounded by
	\begin{equation}\label{eq:ExpDec}
	|f(x)| \leq L \e^{-\alpha \left| x \right|}, \quad \forall x \in \R,  
	\end{equation}
	with some  $\alpha,L$ > 0.
	Then the error of $2N+1$ term sinc interpolation of $f(x)$ by $C_N\{f, h\} (x)$, 
	satisfies the following estimate
	\begin{equation}\label{eq:SAErrEst}
	\begin{split}
	\sup\limits_{x \in \R} \left |f(x) - C_N \{f, h\}(x)\right| \leq c \Eg_N, \\
	\Eg_N = N^{1/2}\e^{-\sqrt{\pi d \alpha N}},
	\end{split}
	\end{equation}
	provided that 
	\begin{equation}\label{eq:SAh}
	h = \sqrt{\frac{\pi d}{\alpha N}}.
	\end{equation}
	Here $c>0$ is some constant dependent on $f, d, \alpha$ and independent on $N$.
}
In this paper we extend the results of the above theorem to a class of algebraically decaying functions on $\R$.
All theoretical considerations are given in sections 1,2. Section 3 is devoted to numerical examples and discussion. 

\section{Interpolation of functions with algebraic decay on real line}
In this section we study the convergence of sinc interpolation for the class of algebraically decaying functions. Specifically, we consider the situation when the  function  $f(x)$ satisfies 
\begin{equation}\label{eq:AlgDec}
|f(x)| \leq  \frac{L}{1 +  |x|^\alpha}, \quad \forall x \in \R
\end{equation}
instead of inequality \eqref{eq:ExpDec}, convenient for the classical sinc methods \cite{Stenger:1993:NMB}.

\begin{theorem}\label{thm:SA_alg_decay}
	Assume that the function $f \in H^1(D_d)$ has an algebraic decay defined by \eqref{eq:AlgDec} with some $\alpha > 1$, $L > 0$. 
	Then the error of $2N+1$-term sinc interpolation \eqref{eq:CN}
	satisfies the following estimate 
	\begin{equation}\label{eq:err_est}
	\begin{split}
	&\sup\limits_{x \in \R} \left |f(x) - C_N \{f, h\}(x)\right| \leq c\Eg_N, \quad \forall x \in \R, \\
	&\Eg_N =\frac{\alpha^{\alpha}(N+1)^{1-\alpha}}{(\alpha-1)(\pi d)^\alpha}
	\left(
	\bW \left(
	\frac{\pi d }{\alpha}
	\left(\frac{\alpha - 1}{\pi d}\right)^{\frac{1}{\alpha}} 
	(N+1)^{\frac{\alpha - 1}{\alpha}}
	\right)
	\right)^{\alpha},
	\end{split}
	\end{equation}
	provided that $h$ in \eqref{eq:CN} is chosen as 
	\begin{equation}\label{eq:h_exact}
	h=\frac{\pi d }{\alpha}
	\left(
	\bW \left(
	\frac{\pi d }{\alpha}
	\left(\frac{\alpha - 1}{\pi d}\right)^{\frac{1}{\alpha}} 
	(N+1)^{\frac{\alpha - 1}{\alpha}}
	\right)
	\right)^{-1}.
	\end{equation}
	Here $\bW[\cdot]$ denotes a positive branch of the Lambert-W function, 
	$c =   c_1 N_1(f,D_d) + 2 L$  and $c_1 > 1$ is the constant 
	independent of $N$:
	\begin{equation}\label{eq:c_1}
	c_1 =   
	\dfrac
	{
		(\pi d)^{2(\alpha-1)}(\alpha-1)^{2}
	} 
	{
		(\pi d)^{2(\alpha-1)}(\alpha-1)^{2}
		-
		\alpha^{2\alpha}
		\bW^{  2 \alpha } \left(
		\frac{\pi d }{\alpha}
		\sqrt[\alpha]{\frac{\alpha - 1}{\pi d}}
		\right)
	}.
	\end{equation}
\end{theorem}
\textit{Proof.}
	For any fixed $h$ the error of sinc interpolation can be represented as follows \cite[equation (3.1.29)]{Stenger:1993:NMB}
	$$
	|f(x) - C_N \{f, h\}(x)| \leq |f(x) - C_\infty \{f, h\}(x)| + \suml_{|k|>N} \left |f(kh) \right |
	$$
	Bound of the first term on the right-hand side 
	of this formula  was obtained in Theorem 3.1.3 from \cite{Stenger:1993:NMB}. For $x \in \R$ this term satisfies
	\begin{equation}\label{eq:disc_err}
	|f(x) - C_\infty \{f, h\}(x)| \leq 
	\frac{N_1(f,D_d) }{2\pi d \sinh{\frac{\pi d}{h}}} \leq
	\frac{c_1 N_1(f,D_d) }{\pi d} \e^{-\frac{\pi d}{h}},
	\end{equation}
	where $c_1>1$ is some constant to be determined later.
	For the second term we get 
	\begin{equation}\label{eq:trunc_err}
	\begin{aligned}
	\suml_{|k|>N} \left |f(kh) \right | 
	&\leq 2 L  \suml_{k=N+1}^\infty (kh)^{-\alpha} 
	\leq 2L \intl_{N+1}^\infty (th)^{-\alpha} d t\\
	&\leq \frac{2L (N+1)^{1-\alpha}}{(\alpha - 1 )h^\alpha}.
	\end{aligned}
	\end{equation}
	The above sequence of inequalities is justified as long as $f(x)$   satisfy \eqref{eq:AlgDec} with some $\alpha>1$. 
	For such $f(x)$, truncation error \eqref{eq:trunc_err} decays algebraically as $N \rightarrow \infty$. In order to balance it with  exponentially decaying discretization error  \eqref{eq:disc_err} one needs to solve for $h$ the equation
	\begin{equation}\label{eq:speed_balance_eq}
	\frac{\e^{-\frac{\pi d}{h}}}{c_2}=  \frac{ (N+1)^{1-\alpha}}{(\alpha - 1 )h^\alpha}. 
	\end{equation}
	
	Let $s = \frac{\pi d }{\alpha} h^{-1}$ and assume that $c_2>0$ is some fixed parameter. Then, equation \eqref{eq:speed_balance_eq} takes the form 
	\[
	\frac{\pi d }{\alpha}   
	\left( \frac{\alpha - 1}{c_2}
	(N+1)^{\alpha-1}
	\right)^{\frac{1}{\alpha}}
	=
	s \e^s,
	\]
	which has a unique solution
	\[
	s = \bW \left(\frac{\pi d }{\alpha}
	\left( 
	{\frac{\alpha - 1}{c_2} (N+1)^{\alpha-1}}
	\right)^{\frac{1}{\alpha}}
	\right). 
	\]
	Next, we set $c_2 = \pi d$ 
	and substitute back the expression for $s$ in terms of $h$ to obtain \eqref{eq:h_exact}.
	The proof of \eqref{eq:err_est} is straightforward
	\begin{align*}
	|f(x) &- C_N \{f, h\}(x)| \leq 
	\left( 
	c_1 N_1(f,D_d)  + 2 L
	\right)
	\frac{ (N+1)^{1-\alpha}}{(\alpha - 1 )h^\alpha} \\
	&\leq 
	c 
	\frac{\alpha^{\alpha}(N+1)^{1-\alpha}}{(\alpha-1)(\pi d)^\alpha}
	\left(
	\bW \left(
	\frac{\pi d }{\alpha}
	\left(\frac{\alpha - 1}{\pi d}\right)^{\frac{1}{\alpha}} 
	(N+1)^{\frac{\alpha - 1}{\alpha}}
	\right)
	\right)^{\alpha} .
	\end{align*}
	Now, let us come back to the determination  of $c_1$. 
	The smallest $c_1$ suitable for \eqref{eq:disc_err} can be defined as follows 
	\[
	c_1  = 
	\sup_{N \in \Z_{+}}
	\left\{ 
	\frac{e^{\frac{\pi d}{h}}}{2 \sinh{\frac{\pi d}{h}}} 
	\right\}  
	= 
	\max_{N \in \Z_{+}}
	\left( 
	{1-e^{-\frac{2\pi d}{h}}} 
	\right)^{-1}. 
	\]
	Its not hard to see that the maximum is attained at $N=0$. 
	Therefore, the value of $c_1$: 
	\[
	c_1  = \left( 1 - \exp\left(- 2 \alpha 
	\bW \left(
	\frac{\pi d }{\alpha}
	\sqrt[\alpha]{\frac{\alpha - 1}{\pi d}}
	\right)
	\right) \right)^{-1} 
	\]
	is clearly greater than one, for any $\alpha > 1$, $d>0$. 
	%
	To get \eqref{eq:c_1} we apply the identity $\exp{\left(-\bW(x)\right) } = \bW(x)/x$ to the above formula for $c_1$ and rearrange the result accordingly 
	\begin{align*}
	c_1  
	& = \left( 
	1 - 
	\frac{\alpha^{2\alpha}}{(\pi d)^{2(\alpha-1)}(\alpha-1)^{2} }
	\left(
	\bW \left(
	\frac{\pi d }{\alpha}
	\sqrt[\alpha]{\frac{\alpha - 1}{\pi d}}
	\right)
	\right)^{ 2 \alpha }
	\right)^{-1} \\
	& =  
	\dfrac
	{
		(\pi d)^{2(\alpha-1)}(\alpha-1)^{2}
	} 
	{
		(\pi d)^{2(\alpha-1)}(\alpha-1)^{2}
		-
		\alpha^{2\alpha}
		\bW^{  2 \alpha } \left(
		\frac{\pi d }{\alpha}
		\sqrt[\alpha]{\frac{\alpha - 1}{\pi d}}
		\right)
	}.
	\end{align*}
{\ \hfill$\Box$}

The presence of $\bW(x)$ in estimate \eqref{eq:err_est} makes it harder to perceive the asymptotic behavior of the interpolation error intuitively.  
To fix that we recall a well-established result \cite{Corless1996} on the asymptotic properties of $\bW(x)$, valid for any $x > e$:  
\[
\ln{x} - \ln{(\ln{x})} + \frac{\ln{(\ln{x})}}{2 \ln{x}} 
\leq \bW(x) \leq
\ln{x} - \ln{(\ln{x})} + \frac{e\ln{(\ln{x})}}{(e-1) \ln{x}}.
\]
By using the above inequality along with the definition of $\bW(x)$ and \eqref{eq:speed_balance_eq} we transform \eqref{eq:err_est} in the following way 
\begin{multline*}
|f(x) - C_N \{f, h\}(x)| 
\leq \frac{c}{\e^{\alpha s}} 
\leq c
\left(
\frac{
	\ln\left(
	\frac{\pi d }{\alpha}
	\left(\frac{\alpha - 1}{\pi d}\right)^{\frac{1}{\alpha}} 
	(N+1)^{\frac{\alpha - 1}{\alpha}}
	\right) 
}
{	\frac{\pi d }{\alpha}
	\left(\frac{\alpha - 1}{\pi d}\right)^{\frac{1}{\alpha}} 
	(N+1)^{\frac{\alpha - 1}{\alpha}}
}
\right)^\alpha  
\\
\leq 
\frac{c}{\left(\pi d\right )^{\alpha - 1}}
\left(
\frac{N + 1}{\alpha - 1 }
\right )^{1-\alpha}
\ln^{\alpha}
\left(
\pi d
\left(\frac{\alpha - 1}{\alpha^\alpha}\right)^{\frac{1}{\alpha-1}}
(N+1)
\right);
\end{multline*}
whence it is clear that the error of sinc interpolation provided by Theorem \ref{thm:SA_alg_decay} is asymptotically equal to $(N+1)^{1-\alpha} \ln^\alpha(N+1)$ as $N \rightarrow \infty$.
To analyze the error 
for small $N$ we note that, in the view of \eqref{eq:speed_balance_eq},  $\Eg_N$ is bounded by the exponent with a strictly decreasing negative argument. Consequently,  for any $\alpha>1$, $x \in \R$,  the error $\sup_{x \in \R}|f(x) - C_N \{f, h\}(x)|$ lies within the interval $[0, c]$ and decreases as $N\rightarrow \infty$.

One might conclude from the foregoing analysis that a simple asymptotic formula $W(x) \approx \ln(x)$ can be used to redefine $h$  \eqref{eq:h_exact} in terms of logarithms, which are computationally more favorable than the Lambert-W function.
To explore this possibility we set 
\[\begin{aligned}
h=&\frac{\pi d }{\alpha}\left(
\ln \left(
\frac{\pi d }{\alpha}
\left(\frac{\alpha - 1}{c_2}\right)^{\frac{1}{\alpha}} 
(N+1)^{\frac{\alpha - 1}{\alpha}}
\right)
\right)^{-1} ,
\end{aligned}\]
and study the corresponding error terms of the approximation. 
Discretization error \eqref{eq:disc_err} 
is positive and monotonically decreasing in $N$ for any $c_2>0$, since $h$ is monotonic. 
The 
principal part $\frac{ (N+1)^{1-\alpha}}{(\alpha - 1 )h^\alpha}$ of truncation error \eqref{eq:trunc_err} 
has one global maximum at $N=N_0$:
\[
N_0=
\left( {\frac {\alpha}{\pi\,d}} \right) ^{{\frac {\alpha}{\alpha-1}}} 
{\exp\left ({{\frac {\alpha}{\alpha-1}}}\right )}
{ \left( \frac{\alpha - 1}{c_2} \right) }^{-\frac{1}{\left( \alpha-1 \right)}}
-1.
\]
To guarantee a monotonous decrease of the truncation error for all $N \geq 0$ we must require $N_0 = 0$, which yields $c_2 = {\left( \alpha-1\right) \left( {\frac{\pi d} {\alpha {\rm e}}} \right)^{\alpha}}.$ 
The aforementioned formula for $h$ is thereby reduced to
\begin{equation}\label{eq:h_ln}
h = \frac{\pi d}{\alpha + (\alpha - 1)\ln{\left( N + 1\right)}}.
\end{equation}
For such $h$, the error of sinc interpolation will be bounded by \eqref{eq:err_est} with
\begin{equation}\label{eq:E_ln}
\Eg_N = \frac{ (N+1)^{1-\alpha}}{(\alpha - 1)(\pi d)^\alpha }
\left(
\alpha + (\alpha - 1)\ln{\left( N + 1\right)}
\right)^{\alpha},	
\end{equation}
and $c=(\alpha - 1)\left( \frac{\pi d}{\alpha e}\right)^\alpha N_1(f,D_d) +  2L$. The main concern with \eqref{eq:E_ln}, is the presence of additional summand $\alpha$ when compared to \eqref{eq:err_est}. 
\begin{remark}
	The definition of $h$ from Theorem \ref{thm:SA_alg_decay} can not be simplified by adopting $W(x) \approx \ln(x)$, since such simplification, as described by \eqref{eq:h_ln}, \eqref{eq:E_ln}, would make the approximation method ineffective for large $\alpha$. 
\end{remark}


With an additional a priori knowledge about $f(x)$ we should be able to improve the convergence properties of $C_N \{f, h\}(x)$ described by Theorem \ref{thm:SA_alg_decay}. 
The following improvement of \eqref{eq:err_est} offers a more realistic balance of discretization and truncation errors, presuming that both $N_1(f,D_d)$ and $L$ are known.  
\begin{corollary}\label{thm:err_est_h_exact_optimized}
	Assume that the function $f(x)$ satisfies the conditions of Theorem \ref{thm:SA_alg_decay}.  If 
	\begin{equation}\label{eq:h_exact_optimized}
	h=\frac{\pi d }{\alpha}\left(
	\bW \left(
	\frac{\pi d }{\alpha}
	\left(\frac{N_1(f,D_d)(\alpha - 1)}{\pi d L}\right)^{\frac{1}{\alpha}} 
	(N+1)^{\frac{\alpha - 1}{\alpha}}
	\right)
	\right)^{-1},
	\end{equation}
	then the error of sinc interpolation fulfils estimate \eqref{eq:err_est}, with $c=(c_1 + 2)L$ and $\Eg_N$ given by 
	\[
	\Eg_N =\frac{(N+1)^{1-\alpha}}{(\alpha-1)}
	h^{-\alpha}.	
	\]
\end{corollary}
\noindent Formula \eqref{eq:h_exact_optimized} was obtained in the same way as \eqref{eq:h_exact}, except this time we set 
\[
c_2 = \frac{\pi d L}{N_1(f,D_d)}.
\]

\section{Interpolation of functions with algebraic decay in the strip}

Corollary \ref{thm:err_est_h_exact_optimized} is difficult to apply as it is, because the evaluation of $N_1(f,D_d)$ requires computation of the contour integral over $\partial D_d$. 
In order to make this result more applicable we note, that if $f \in H^1(D_r)$, for some $r>0$, then $\lim\limits_{x\rightarrow\pm\infty}f(x + i y) =0$
uniformly with respect to $y \in [d,d]$, for all $d \in (0, r)$  \cite[Proposition 6]{Butzer2013}. 
Hence, for any $r>0$ there exist a nonempty subspace of $H^1(D_r)$, such that its  elements $f$ satisfy 
\begin{equation}\label{eq:Alg_dec_strip}
|f(z)| \leq  \frac{L}{1 +  |z|^\alpha}, \quad \forall z \in D_d,
\end{equation}
with some $d \in (0, r)$. 
\begin{theorem}\label{thm:SA_alg_decay_strip}
	Assume that the function $f(z)$ is analytic in the horizontal strip $D_d$, $d>0$. 
	If $f(z)$ is bounded  by \eqref{eq:Alg_dec_strip}  with some $\alpha > 1$, $L > 0$,
	then the error of sinc interpolation \eqref{eq:CN}
	satisfies the following estimate 
	\begin{equation}\label{eq:err_est_no_N1}
	\begin{split}
	&\sup\limits_{x \in \R} \left |f(x) - C_N \{f, h\}(x)\right| \leq c\Eg_N, \\
	&\Eg_N =\frac{\alpha^{\alpha}(N+1)^{1-\alpha}}{(\alpha-1)(\pi d)^\alpha}
	h^{\alpha},
	\end{split}
	\end{equation}
	provided that 
	\begin{equation}\label{eq:h_no_N1}
	h=
	\frac{\pi d }{\alpha}\left(
	\bW \left(
	\frac{\pi d }{\alpha}
	\left(
	\frac{2 \beta (\alpha - 1)}
	{\pi d}
	\right)^{\frac{1}{\alpha}} 
	(N+1)^{\frac{\alpha - 1}{\alpha}}
	\right)
	\right)^{-1},
	\end{equation}
	with $\beta = \min\left\{
	\tfrac{2}{\sinc{\left ( \alpha^{-1} \right )}},
	\left(\frac{2}{d}\right)^{\alpha - 1} B\left(\frac{\alpha}{2}-\frac{1}{2}, \frac{\alpha}{2}+\frac{1}{2} \right)
	\right\}.
	$
	Here $B(\cdot,\cdot)$ is the beta function, $c =  c_1 L$ and $c_1$ is the constant dependent on $\alpha, d$.
\end{theorem}
\textit{Proof.}
For small values of $d$ we proceed as follows
	\begin{equation}\label{eq:N1_h_b_est}
	\intl_{-\infty}^{+\infty}\left |f(x + i d)\right | dx \leq 
	\intl_{-\infty}^{+\infty} \frac{L dx }{1+|x+id|^\alpha}  =
	2L\intl_{0}^{+\infty} \frac{dx}{1+(x^2+d^2)^{\frac{\alpha}{2}}}  
	\end{equation}
	\[
	\intl_{0}^{+\infty} \frac{dx}{1+(x^2+d^2)^{\alpha/2}}  \leq 
	\intl_{0}^{+\infty} \frac{dx}{1+x^\alpha}  =
	\lim\limits_{x \rightarrow \infty}{\frac {x{\ \Phi} \left( -{x}^{\alpha},1,{\alpha}^{-1} \right) 
		}{\alpha}}
	\]
	\[
	=\lim\limits_{x \rightarrow +\infty}\left |{\frac {x{\ \Phi} \left( -{x}^{\alpha},1,{\alpha}^{-1} \right) 
		}{\alpha}}\right |
	=
	\lim\limits_{\substack{
			\Re{z} \rightarrow +\infty \\
			\Im{z} \rightarrow 0 
	}}
	\left| {\frac {z{\Phi} \left( -{z}^{\alpha},1,{\alpha}^{-1} \right) 
		}{\alpha}} \right|. 
	\]
	Here $\Re{z}$ and $\Im{z}$ is real and imaginary part of $z$ correspondingly. To evaluate the last limit we employ Corollary 1 from \cite{Ferreira2017}.  It offers a convergent expansion of the Hurwitz-Lerch zeta function ${\Phi(z,s,a)}$, when its second parameter $s$ is an integer number 
	\begin{equation}\label{eq:LerchExpansion}
	z{\Phi} \left( {z}^{\alpha},1,\frac{1}{\alpha} \right) = \pi \left( sgn\left\{\Arg(\alpha\ln(z))\right \} i + \cot{\frac{\pi}{\alpha}} \right) - \suml_{k=1}^{\infty}\frac{z^{1-\alpha k}}{1/\alpha - k}.
	\end{equation}
	The expression on the right of \eqref{eq:LerchExpansion} is bounded and uniformly convergent to the left-hand side for any $\alpha>1$, $|z|>1$,  such that $z^\alpha \notin (-\infty,-1)\cup (1,\infty)$. 
	Therefore
	\[
	\lim\limits_{\substack{
			\Re{z} \rightarrow +\infty \\
			\Im{z} \rightarrow 0 
	}}
	\left| {\frac {z{\Phi} \left( -{z}^{\alpha},1,{\alpha}^{-1} \right) 
		}{\alpha}} \right|  = \frac{\pi}{\alpha} \sqrt{1+\cot^2{\frac{\pi}{\alpha}} } - 
	\frac{1}{\alpha}\suml_{k=1}^{\infty} 
	\lim\limits_{\substack{
			\Re{z} \rightarrow +\infty \\
			\Im{z} \rightarrow 0 
	}}
	\frac{z^{1-\alpha k}}{1/\alpha - k},
	\]
	which leads us to the bound
	\begin{equation}\label{eq:N1_est_small_d}
	\intl_{-\infty}^{+\infty}\left |f(x + i d)\right | dx \leq  \frac{2\pi L}{\alpha}\sqrt{1+\cot^2{\frac{\pi}{\alpha}} } = 
	2L\sinc^{-1}{\left ( \frac{1}{\alpha} \right )}.
	\end{equation}
	For large $d$, the integral from \eqref{eq:N1_h_b_est} can be estimated as follows 
	\begin{align*}
	\intl_{0}^{+\infty} \frac{1}{1+(x^2+d^2)^{\alpha/2}} dx \leq & 
	\intl_{0}^{+\infty} \frac{1}{(x^2+d^2)^{\alpha/2}} dx \\
	&=\frac {\sqrt {\pi} d^{1-\alpha}\Gamma \left( \left (\alpha-1\right )/2 \right) }
	{2 \Gamma	\left(  \alpha/2 \right) } \\
	&= \frac {d^{1-\alpha}\Gamma \left( \left (\alpha-1\right )/2 \right) \Gamma \left( \left (\alpha+1\right )/2 \right)}{2^{2-\alpha}\Gamma(\alpha)}\\
	\leq &
	\frac{1}{2}
	B \left(\frac{\alpha}{2}-\frac{1}{2}, \frac{\alpha}{2}+\frac{1}{2} 
	\right)
	\left(\frac{2}{d}\right)^{\alpha - 1}. 
	\end{align*}
	To obtain the above estimate we used a well-known multiplication theorem \cite[p. 4]{Bateman1953_1} for Gamma function $\Gamma(\cdot)$.
	The next bound is a direct consequence of the above formula and \eqref{eq:N1_h_b_est}
	\begin{equation}\label{eq:N1_est_large_d}
	\intl_{-\infty}^{+\infty}\left |f(x + i d)\right | dx \leq  
	L B \left(\frac{\alpha}{2}-\frac{1}{2}, \frac{\alpha}{2}+\frac{1}{2} \right)
	\left(\frac{2}{d}\right)^{\alpha - 1} \hspace*{-6mm}
	.
	\end{equation}
	
	By combining bounds \eqref{eq:N1_est_small_d}, \eqref{eq:N1_est_large_d} and taking in to account the fact that the expression on the right of \eqref{eq:Alg_dec_strip} is invariant with respect to  $z \rightarrow \bar{z}$ we arrive at the following estimate 
	\begin{equation*}
	N_1(f,D_d) \leq 
	2 L
	\min\left \{
	\dfrac{2}{\sinc{\left ( \frac{1}{\alpha} \right )}},
	\left(\frac{2}{d}\right)^{\alpha - 1}\hspace*{-5mm} B\left(\frac{\alpha}{2}-\frac{1}{2}, \frac{\alpha}{2}+\frac{1}{2} \right)
	\right \}.
	\end{equation*}
	To finalize the proof, we evaluate \eqref{eq:h_exact_optimized}  assuming that the value of $N_1(f,D_d)$ is equal to its estimate provided by the previous formula. 
	This will get us \eqref{eq:h_no_N1}.
{\ \hfill$\Box$}

\section{Examples and discussion}
In this section we consider several examples of the developed approximation method. 
As a measure of experimental error we use a discrete norm
\[
err = \max\limits_{\forall x \in X} \left| f(x) - C_N \{f, h\}(x)\right| , 
\]
defined on a uniform grid $X = \left \{{j h}/{2} \ \middle|\  j= \overline{-2N,2N} \right \}$. 
With such choice of $X$ the specified discrete norm ought to capture the contribution from both the descretization and truncation parts of the error. To experimentally check the convergence of $C_N \{f, h\}(x)$ we repeat the approximation procedure on a sequence of grids determined by 
\[
N_i \in \{1,2,4,8,16,32,64,128,256,512,1024\},
\] 
and the corresponding $h_i$ evaluated by one of the formulas \eqref{eq:h_exact}, \eqref{eq:h_exact_optimized} or \eqref{eq:h_no_N1}.

\begin{example}\label{ex:1}
	Let	
	\[
	f(x) = \dfrac{4}{2+x^{2 a}},
	\]
	where $a \geq 2$ is integer.
	Then, the largest possible value of $d$ such that $f(x)$ remains analytic in $D_d$, is equal to 
	$\sqrt[2\alpha]{2} \sin{\frac{\pi}{2a}}$. For the purpose of the illustration we set  $d=0.9\sqrt[2a]{2} \sin{\frac{\pi}{2a}}$, $a=2$, then $N_1(f,D_d) \approx 17.05467564$, $L \approx 4$, $\alpha = 4$.  The behaviour of an error $err(x) = f(x) - C_{32} \{f, h\}(x)$ for the values of $h$, calculated by three different formulas \eqref{eq:h_exact}, \eqref{eq:h_exact_optimized}, \eqref{eq:h_no_N1}, is depicted in Fig. \ref{fig:ex1}.
	\begin{figure}[ht]
		\begin{center}
			\includegraphics[width=0.32\linewidth]{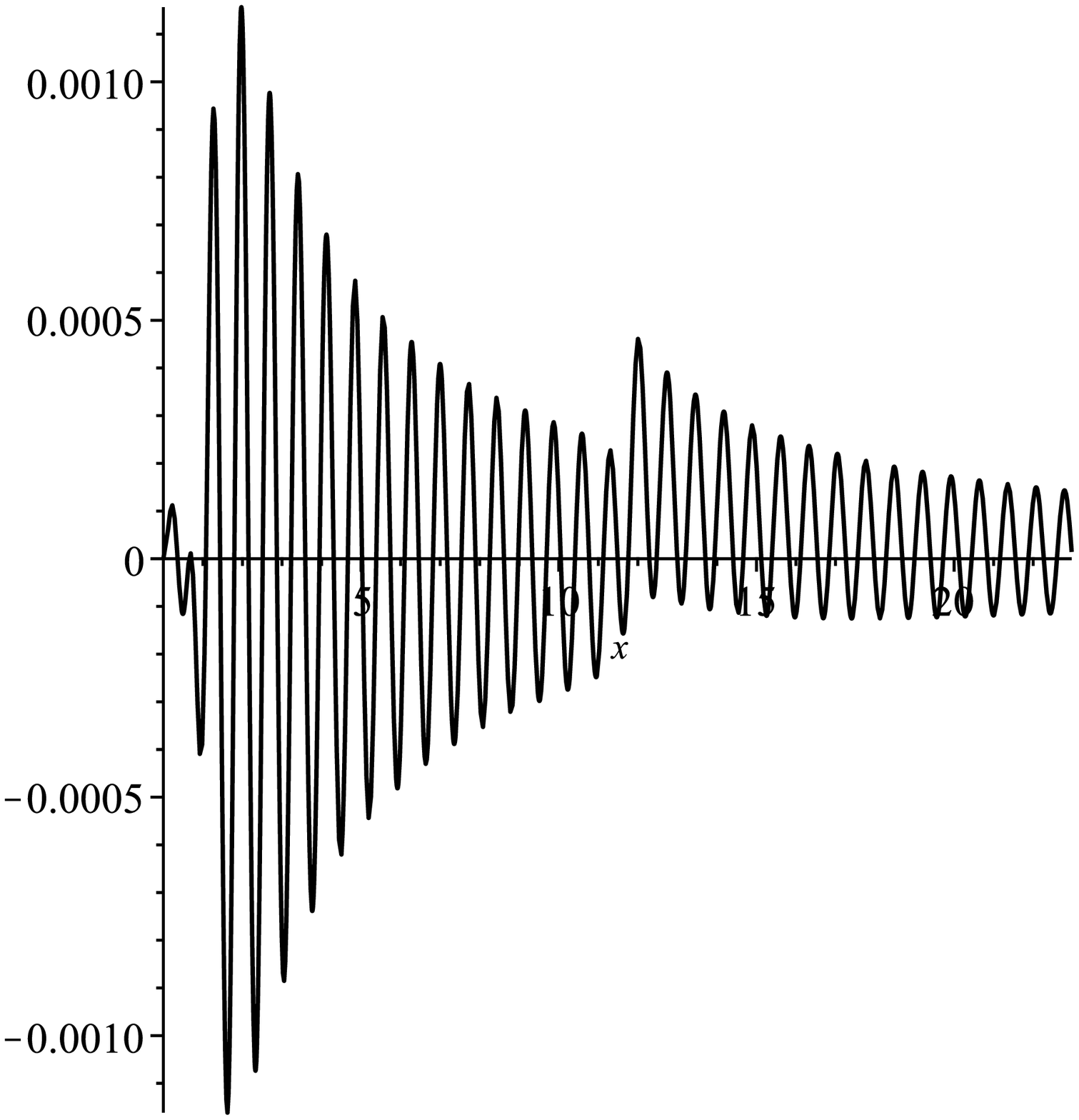}
			\includegraphics[width=0.32\linewidth]{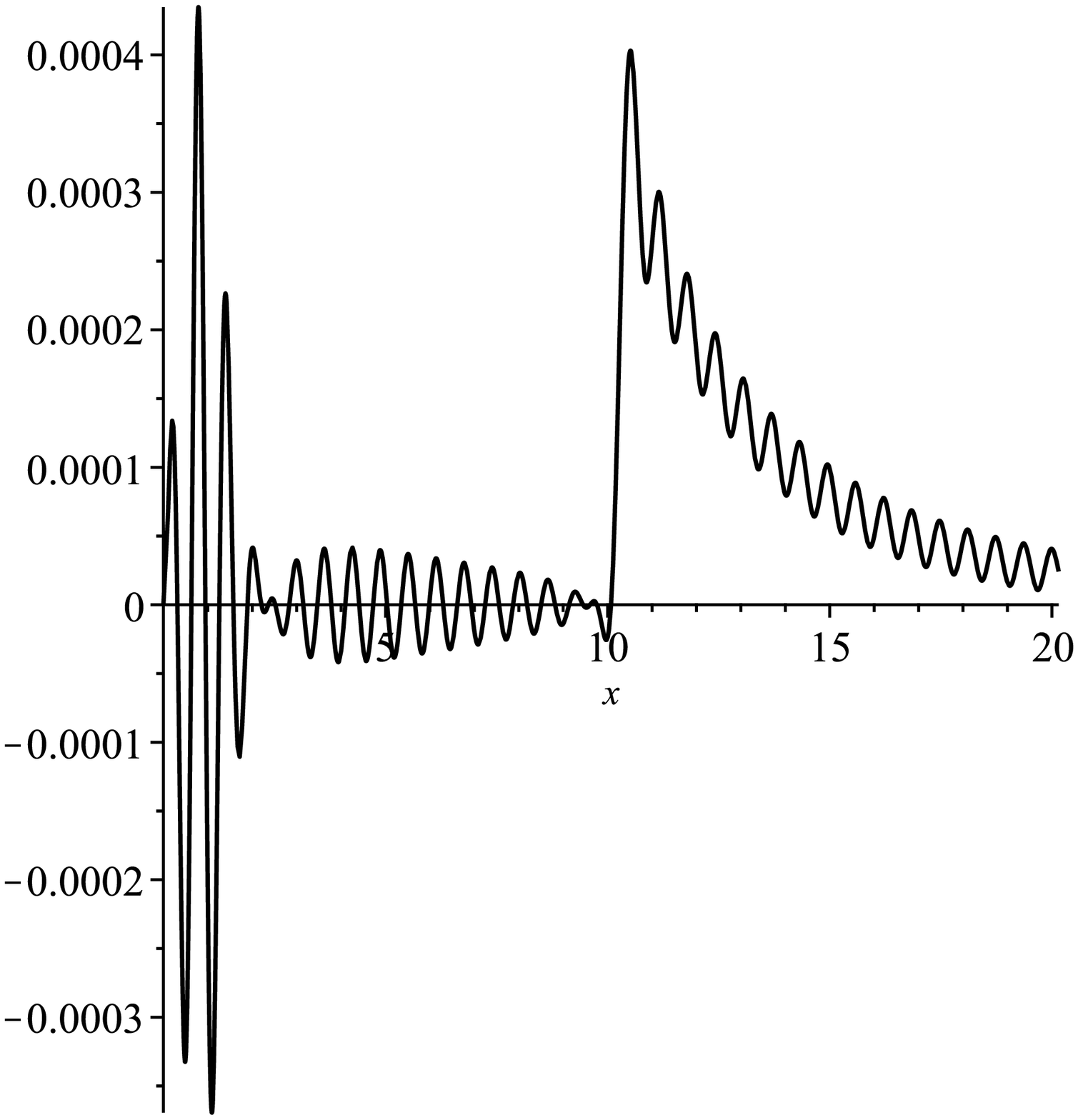}
			\includegraphics[width=0.32\linewidth]{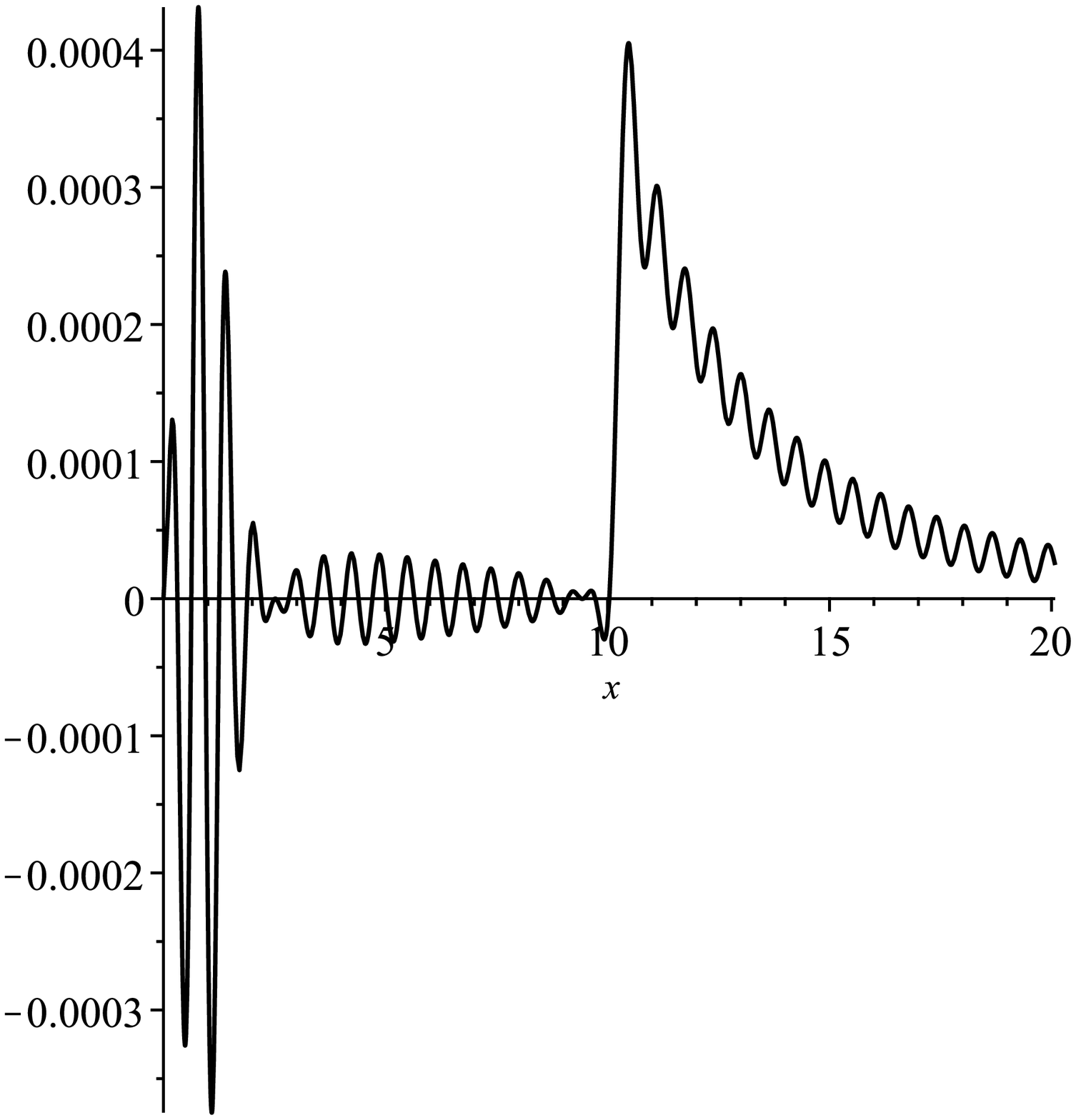}
		\end{center}
		\caption{Graphs of $err(x) = f(x) - C_{32} \{f, h\}(x)$ from Example \ref{ex:1} for $h$ calculated by \eqref{eq:h_exact} -- left graph, \eqref{eq:h_exact_optimized} -- central graph and \eqref{eq:h_no_N1} -- graph on the right. }
		\label{fig:ex1}
	\end{figure}
	Predictably, the maximum of $err(x)$ for $h=0.3149022805$ calculated by \eqref{eq:h_exact_optimized} (see central plot from Fig. \ref{fig:ex1}) is superior to the error with  $h=0.3589479879$ calculated by \eqref{eq:h_exact} (left plot from Fig. \ref{fig:ex1} ). The value of $h$ calculated by \eqref{eq:h_no_N1} is close to the one obtained from \eqref{eq:h_exact_optimized}, that is why the error function $err(x)$ (see plot on the right from Fig. \ref{fig:ex1}) is close to $err(x)$ obtained with help of \eqref{eq:h_exact_optimized}.
	One can see a discernible  spike in the error function from central plot of Fig. \ref{fig:ex1} at  $x_0 = N_6 h \approx 10.0769$. 
	The values of $err(x)$ on the left of $x_0$ corresponds to the discretization error, whilst  the values on the right of $x_0$ corresponds to the truncation error. 
	The magnitude of those errors almost match. 
	This highlight the fact that the chosen $h$ is quite close to the theoretically optimal value.
\end{example}
\begin{example}\label{ex:2}
	In this example we set $f(x) \in H^1(D_d)$ as   
	\[
	f(x) = \dfrac{6 \cos{2 x}}{\left (5 + \cos^2{x}\right)\left(1+x^4\right)},
	\]
	and choose formula \eqref{eq:h_exact} for the evaluation of $h$.
	The function $f(x)$ is meromorphic and bounded in $D_d$ for any $d$ smaller than the imaginary part of zeros of $\left (5 + \cos^2{x}\right)\left(1+x^4\right)$. 
	The zeros of the polynomial part of this expression lie closer to the real line than any zero of $5 + \cos^2{x}$, so $d \le \Im{\sqrt[4]{-1}} = \frac{\sqrt{2}}{2} \approx .707106781186550$. 
	Therefore it is safe to set $d = 0.7$. 
	For given $f(x)$ we can also explicitly find  the parameters of algebraic decay bound \eqref{eq:AlgDec}: $L=f(0)=1$, $\alpha = 4$. 
	
	Note, that for a more general function $f(x)$ the corresponding $L,\alpha$ can be calculated numerically from a sequence of its values. 
	For explicitly given $f(x)$ the possible values of $d$ can be calculated numerically as well, for  example  using \verb|Analytic| routine from Maple \cite{maple_2016}.
	
	The graphs of the approximated function $f(x)$ and the error of its interpolation by $C_{32} \{f, h\}(x)$ are given in Fig. \ref{fig:ex2}. 
	\begin{figure}[ht]
		\begin{center}
			\includegraphics[width=0.49\linewidth]{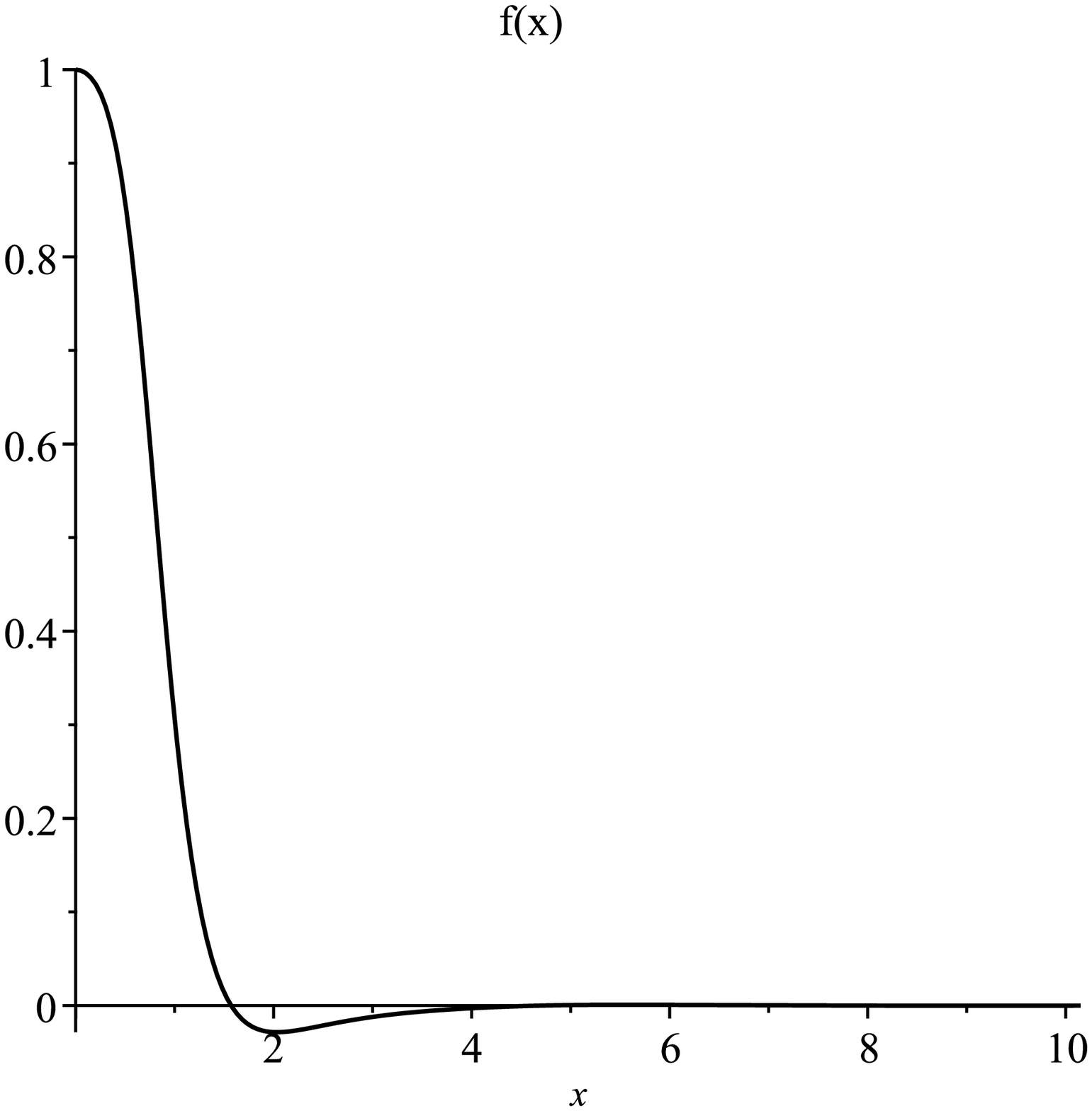}
			\includegraphics[width=0.49\linewidth]{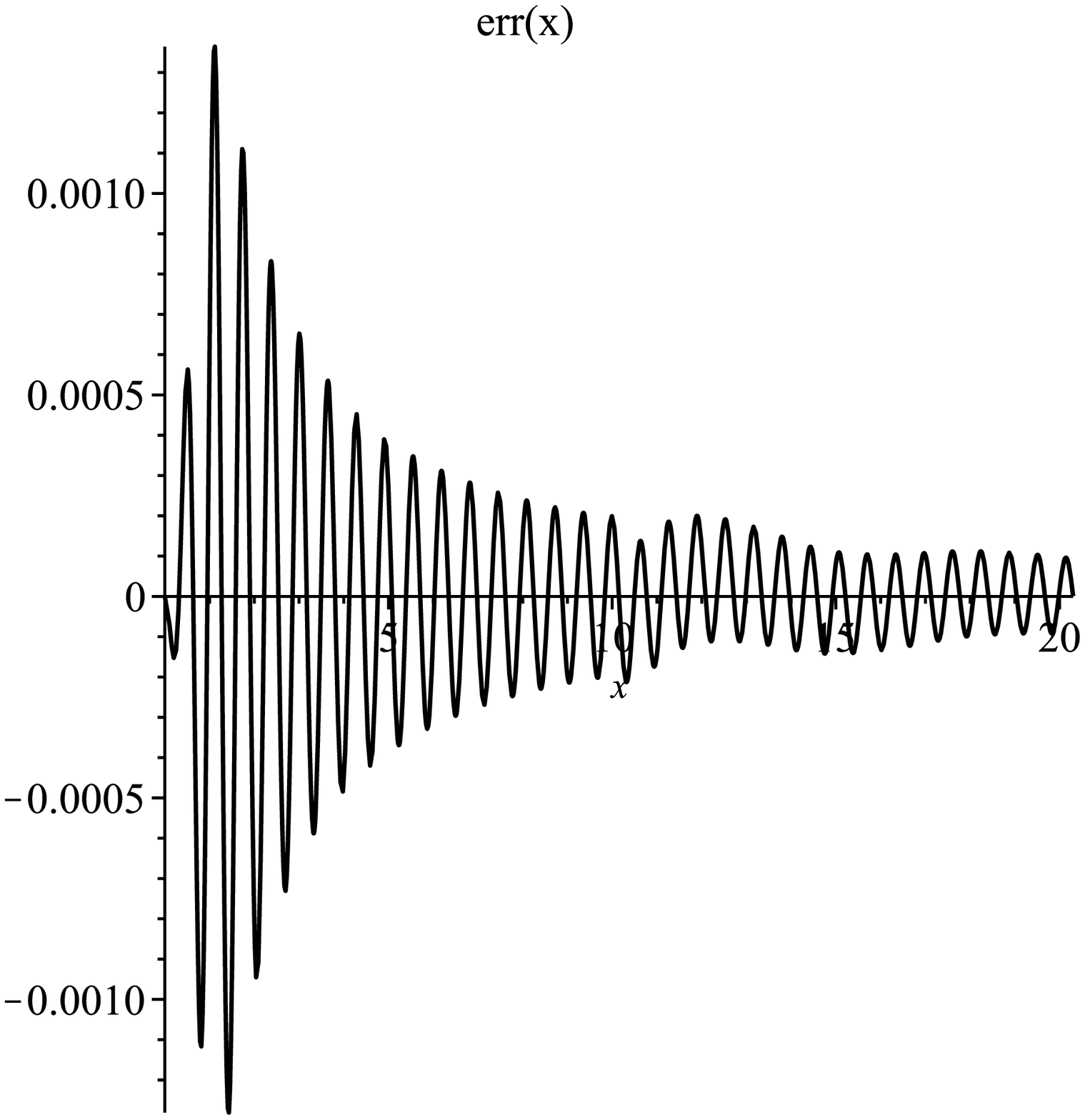}
		\end{center}
		\caption{Graphs of $f(x)$ and $err(x) = f(x) - C_{32} \{f, h\}(x)$ from Example \ref{ex:2}. }
		\label{fig:ex2}
	\end{figure}
	The precise values of $err_i$ for $i = 1,\ldots,11$ are presented in Table \ref{tab:ex2}. 
	Here we additionally supply the theoretical estimate $\Eg_{N_i}$ defined in Theorem \ref{thm:SA_alg_decay} and the value of $c_i = err_i/\Eg_{N_i}$. 
	\begin{table}[ht]
		\caption{Result of the numerical experiments for $f(x)$ from Example \ref{ex:2}. The step size $h$ is calculated by  \eqref{eq:h_exact}, the quantities $\Eg_{N}$ and $c$ are evaluated with help of \eqref{eq:err_est}.}		
		\begin{center}
			\begin{tabular}{|c|c|c|c|c|}
				\hline
				$i$ & $N_i$ & $err_i$ & $\Eg_{N_i}$ & $c_i$ \\
				\hline
1& 2& 6.373770E-02& 3.641222E-02& 1.750448\\
2& 4& 4.011175E-02& 1.904281E-02& 2.106399\\
3& 8& 1.019463E-02& 8.186076E-03& 1.245362\\
4& 16& 3.765622E-03& 2.948999E-03& 1.276915\\
5& 32& 1.368552E-03& 9.160491E-04& 1.493972\\
6& 64& 1.777309E-04& 2.523604E-04& 0.704274\\
7& 128& 7.216260E-05& 6.312895E-05& 1.143098\\
8& 256& 7.698800E-06& 1.460731E-05& 0.527051\\
9& 512& 2.505400E-06& 3.171023E-06& 0.790092\\
10& 1024& 3.281000E-07& 6.528835E-07& 0.502540\\
\hline
			\end{tabular}	
		\end{center}
	\label{tab:ex2}
	\end{table}

\end{example}
	The data from  in Table \ref{tab:ex2} demonstrates that the approximation method presented by Theorem \ref{thm:SA_alg_decay} converges to $f(x)$. 
The magnitude of the observed approximation errors are consistent with the estimate provided by \eqref{eq:err_est}. 
Moreover the estimated value of $c$ from \eqref{eq:err_est} remains bounded by $2.1$ for all $i=\overline{1,6}$. All this prove the effectiveness of the developed method.

%
\begin{small}
\begin{flushleft}
\textsc{Institute of mathematics, NAS of Ukraine,
3, Tereschenkivska St., 01004, Kyiv, Ukraine}\\
\textit{E-mail address}: sytnikd@gmail.com, syntik@imath.kiev.ua
\end{flushleft}
\end{small}

\end{document}